# Existence of periodic orbits in three-dimensional piecewise linear systems

[1,2], Songmei Huan,  [1] Xiao-Song Yang

[1] Department of Mathematics,
[2] Department of Electronics and Information Engineering,
Huazhong University of Science and Technology.
Wuhan, 430074, China

**Abstract** Based on the results about the invariant cones appeared in the literature this paper analyses the existence of periodic orbits in three-dimensional continuous piecewise linear homogeneous systems with two zones, and a necessary and sufficient condition for the existence of periodic orbits of such systems is given.

## 1 Introduction

Many dynamical systems that are used to describe physical process are piecewise smooth. For example, dynamics of some physical systems can be characterized by periods of smooth evolutions interrupted by instantaneous events. Extensive studies about the so called piecewise smooth systems can be found in [3], [10] and [11] and references therein.

Piecewise linear systems as one of the simplest kinds of piecewise smooth systems are successfully used to model physical phenomena and engineering problems, and they can display almost the same dynamical behaviors of general nonlinear systems, such as limit cycles, homoclinic and heteroclinic orbits, strange attractors, and so on. A lot of such examples obtained from control and mechanics can be found in the books [4], [5] and [16]. In nonlinear oscillation theory, the low dimensional continuous piecewise linear systems with one or two separation boundaries in $R^2$ and $R^3$ frequently appear in applications (the reader can see [1], [15], [17] and the references therein to get a comprehensive review about this).

Because piecewise smooth systems are non-differentiable and even discontinuous, the methods from the differentiable dynamical systems theory cannot be applied and specific approaches have to be pursued. As for piecewise linear continuous time system, although explicit solutions can be obtained in each region where one linear system is defined, gluing these solutions globally at the discontinuity boundaries to get a whole picture of the piecewise linear system is a hard problem, and the theory for piecewise linear continuous time systems is also in its infancy, see [5] and [12].

For piecewise linear systems, Freire et al. gave a thorough analysis for two-dimensional systems with two zones in [13] and a similar study for symmetrical two-dimensional systems with three zones in [14]. The case for the continuous piecewise systems in $R^3$ with two zones can be rather complex, as shown in [2].

In [8], the authors found several invariant manifolds foliated by periodic orbits for some degenerate continuous piecewise linear systems with two zones in $R^3$, whose linear parts share a pair of imaginary eigenvalues. Motivated by this work, the existence of periodic orbits for general observable continuous three-dimensional homogeneous continuous piecewise linear systems with two zones is studied in this paper. Several results are given with some examples to illustrate the conclusions about the existence of periodic orbits.

For simplicity, we just discuss the existence of periodic orbits intersecting the separation plane $x_1 = 0$ and concentrate our attention on homogeneous continuous piecewise linear systems in $R^3$ with a separation plane that can be written without loss of generality in the form

$$\dot{x} = \begin{cases} A^+ x & x_1 \geq 0 \\ A^- x & x_1 < 0 \end{cases}, \quad x \in R^3 \tag{1}$$

where the dot denotes derivatives respect to the time $t$, $x = (x_1\ x_2\ x_3)^T \in R^3$ and $A^+, A^-$ are $3 \times 3$ real matrices that satisfy the continuity relation, i.e., $A^+, A^-$ share the second and third columns [6,7]. In addition, if system (1) is **observable**, that is, the observability matrix

$$O = \begin{bmatrix} e_1^T \\ e_1^T \cdot A^- \\ e_1^T \cdot (A^-)^2 \end{bmatrix}$$

has full rank, it is possible to perform a similar transformation to put system (1) into the form

$$\dot{x} = B^{\pm} x = \begin{pmatrix} \delta^{\pm} & -1 & 0 \\ m^{\pm} & 0 & -1 \\ d^{\pm} & 0 & 0 \end{pmatrix} \cdot x$$

with $\delta^{\pm}, m^{\pm}$ and $d^{\pm}$ being the coefficients of the characteristic polynomials of matrices $A^{\pm}$ given by

$$p^{\pm}(\lambda) = \lambda^3 - \delta^{\pm}\lambda^2 + m^{\pm}\lambda - d^{\pm},$$

respectively, see [9] for more details. Hence, in the sequel, without loss of generality, we assume in system (1)

$$A^{\pm} = \begin{pmatrix} \delta^{\pm} & -1 & 0 \\ m^{\pm} & 0 & -1 \\ d^{\pm} & 0 & 0 \end{pmatrix}, \tag{2}$$

also replace $x_2, x_3$ by $y, z$ for convenience.

Note that when $A^{\pm}$ are given in form (2), system (1) is observable.

In the following analysis, we focus on the case when $A^{\pm}$ both have a real and a pair of conjugate complex eigenvalues $\lambda^{\pm}, \alpha^{\pm} \pm \beta^{\pm} j$ with $\beta^{\pm} > 0$. By the characteristic polynomials of matrices $A^{\pm}$, we have

$$\begin{cases} \lambda^{\pm} + 2\alpha^{\pm} = \delta^{\pm} \\ 2\lambda^{\pm}\alpha^{\pm} + (\alpha^{\pm})^2 + (\beta^{\pm})^2 = m^{\pm} \\ \lambda^{\pm} \cdot [(\alpha^{\pm})^2 + (\beta^{\pm})^2] = d^{\pm} \end{cases} \quad (3)$$

By means of the flow of system $\dot{x} = A^{-}x$ with $x_1 \leq 0$, the point $p(0, y_0, z_0)$ with $y_0 > 0$ will be transformed into $q(0, y_1, z_1)$ with $y_1 < 0$, so a left Poincaré half map $P^{-}$ can be defined as $(y_1, z_1) = P^{-}(y_0, z_0)$. Analogously, a right Poincaré half map $P^{+}$ can be defined in some region of the half plane $x_1 = 0, y < 0$ as $(y_2, z_2) = P^{+}(y_1, z_1)$ with $y_1 < 0$ and $y_2 > 0$. Then the composite Poincaré map $P = P^{+} \circ P^{-}$ can be defined from a subset of the half plane $x_1 = 0, y > 0$ into itself, see **Fig.1**.

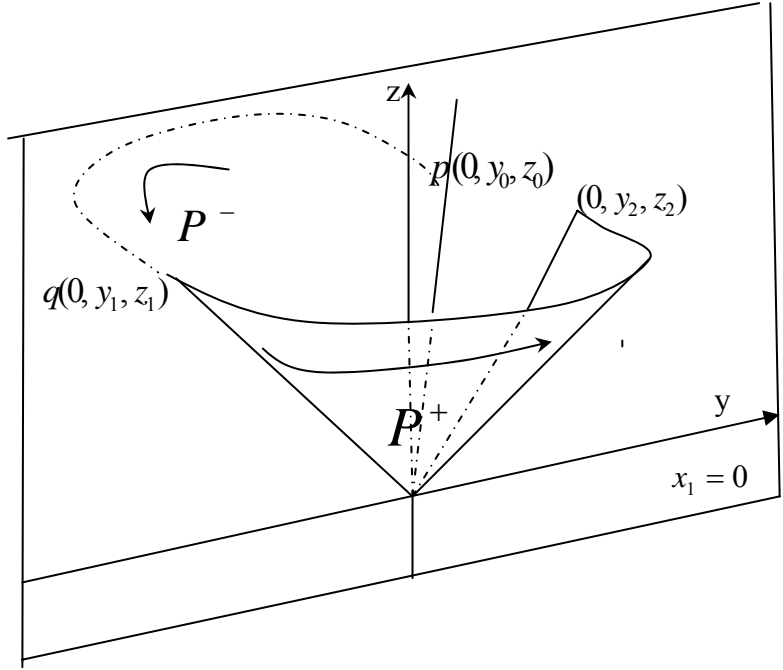

**Fig.1** Poincaré map of system (1)

It is easy to see that if $x(t)$ is a solution of system (1), then $\mu \cdot x(t)$ for $\mu > 0$ is also a solution. So, according to the Poincaré half maps $P^{-}$ and $P^{+}$, we can define two slope transition maps $u_1 = S^{-}(u_0)$ and $v_2 = S^{+}(v_1)$, respectively. Here $u_0$ is the slope of the line $x_1 = 0$, $z = u_0 y$ with $y > 0$, $u_1$ is the image of $u_0$ through $P^-$ with $x_1 = 0$, $z = u_1 y$ where $y < 0$, and $v_2$ (with $x_1 = 0$, $z = v_2 y$ where $y > 0$) is the image of $v_1$ (with $x_1 = 0$, $z = v_1 y$ where $y < 0$) through $P^+$. See **Fig.1** again.

Note that that the composite map $S^+ \circ S^-$ having fixed points, i.e., system (1) having two-zonal invariant cones (see **Definition**) is the necessary condition for the existence of periodic orbits.

**Definition** Let $l$ be a half-line with the origin be its endpoint. If there exists $T > 0$ such that the solution $\Phi(x_0, t)$ of system (1) beginning from any point $x_0 \in l$ at $t = 0$ will return to $l$ again for the first time at $t = T$, i.e. $\Phi(x_0, T) \in l$, then the surface $C$ formed by $l$ rotating along the flow $\Phi(x_0, t)$, $t \in [0, T]$ is called an **invariant cone**. Furthermore, $C$ will be called **one-zonal invariant cone** if it has no intersection with the separation plane $x_1 = 0$, and otherwise be called **two-zonal invariant cone**.

In order to give the statements of the main results, we also introduce the parameters

$$\gamma^+ = \frac{\alpha^+ - \lambda^+}{\beta^+}, \quad \gamma^- = \frac{\alpha^- - \lambda^-}{\beta^-} \tag{4}$$

and the auxiliary function

$$\phi_\gamma(\tau) = 1 - e^{\gamma \tau}(\cos \tau - \gamma \sin \tau)$$

which are crucial for the analysis to be made in the sequel, and will be introduced in details in the next section.

The rest of the paper is organized as follows. In the second section, some preliminaries are provided. In the third section, we revisit a result given by [7] and give a more geometrically intuitive proof. In the fourth section, the statements and proofs of the main results about the existence of periodic orbits are given. Finally, we give two examples with numerical simulations to illustrate the conclusions given in the fourth section.

## 2 Preliminaries

In this section we will establish the Poincaré half-maps in parametric forms for system (1) with $A^\pm$ both having complex eigenvalues. Since the flow of system (1) is made up by matching two linear flows, we can give the Poincaré map by studying each Poincaré half-map separately as a first step, which has been done in [6] and [7]. For the completeness of the paper, here we just gather these conclusions.

First we introduce the properties of the auxiliary function (which is introduced first in [1])

$$\phi_\gamma(\tau) = 1 - e^{\gamma\tau}(\cos\tau - \gamma\sin\tau)$$

as follows

a) $\phi_\gamma(\tau) = \phi_{-\gamma}(-\tau)$ for any $\gamma \in R$ and any $\tau \in R$;

b) There exists a value $\hat{\tau} \in (\pi, 2\pi)$ such that $\phi_{|\gamma|}(\hat{\tau}) = 0$ for any $\gamma \neq 0$;

c) $\phi_\gamma(\tau) = 1 - \cos\tau$, when $\gamma = 0$.

We sketch the graph of the function $\phi_\gamma(\tau)$ in **Fig.2**.

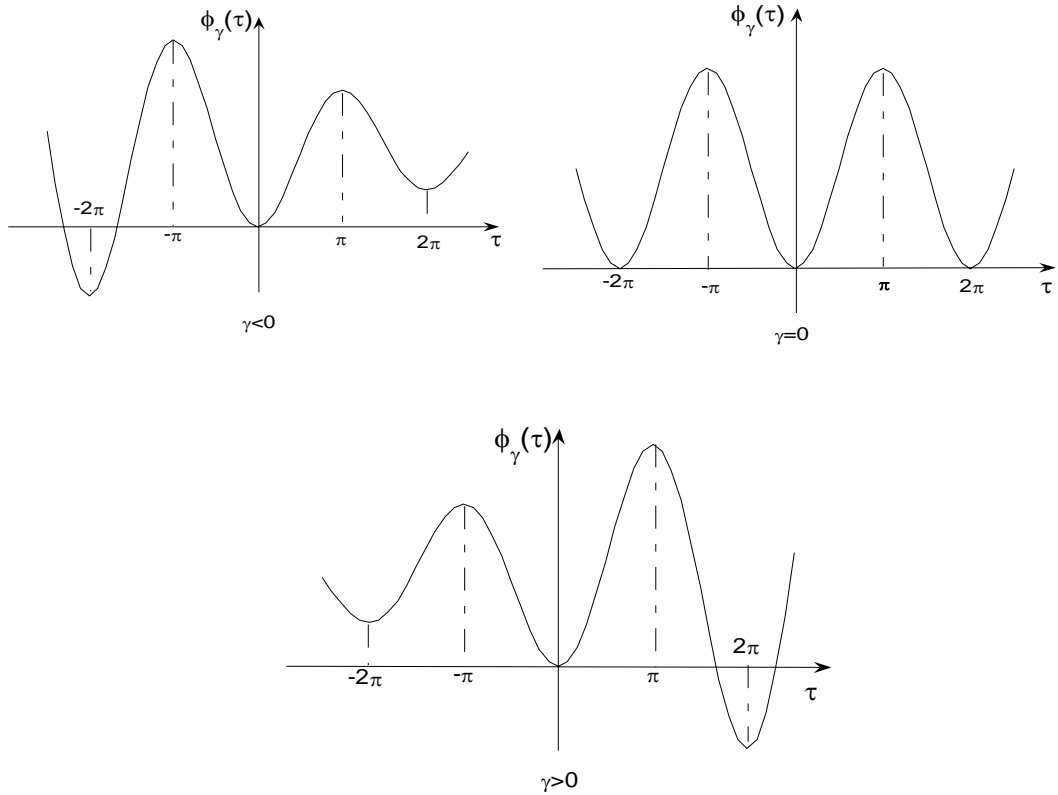

**Fig.2.** The graph of the function $\phi_\gamma(\tau)$.

To computer the Poincaré map of the system $\dot{x} = A^- x$, let us start from the initial point $p(0, y_0, z_0)$ with $y_0 > 0$ at time $t = 0$, then the solution is given by

$$\Phi^-(p,t) = e^{\alpha^- t} \cdot M \cdot \begin{bmatrix} \cos\beta^- t & -\sin\beta^- t & 0 \\ \sin\beta^- t & \cos\beta^- t & 0 \\ 0 & 0 & e^{(\lambda^- - \alpha^-)t} \end{bmatrix} \cdot M^{-1} \cdot \begin{bmatrix} 0 \\ y_0 \\ z_0 \end{bmatrix} \quad (5)$$

where

$$M = \begin{bmatrix} 1 & 0 & 1 \\ \lambda^- + \alpha^- & \beta^- & 2\alpha^- \\ \lambda^-\alpha^- & \lambda^-\beta^- & (\alpha^-)^2 + (\beta^-)^2 \end{bmatrix}.$$

Suppose at time $\tau^-/\beta^-$ the flow returns to the separation plane $x_1 = 0$ for the first time, i.e., $\Phi^-(p, \tau^-/\beta^-) = (0, y_1, z_1)$ with $y_1 < 0$, then by (5) we can get the initial and final coordinate ratios in the parameter forms

$$\begin{aligned} u_0(\tau^-) &= \frac{z_0}{y_0} = \lambda^- + \beta^- \cdot [1 + (\gamma^-)^2] \cdot \frac{e^{\gamma^-\tau^-} \sin \tau^-}{\phi_{\gamma^-}(\tau^-)} \\ u_1(\tau^-) &= \frac{z_1}{y_1} = \lambda^- - \beta^- \cdot [1 + (\gamma^-)^2] \cdot \frac{e^{-\gamma^-\tau^-} \sin \tau^-}{\phi_{-\gamma^-}(\tau^-)} \end{aligned} \quad (6)$$

which not only give the angular behavior of the Poincaré half map $P^-$, but also give the definition of the slope transition map $S^-$ in the parametric form

$$S^- : R \to R, \quad S^-(u_0) = u_1. \tag{7}$$

Furthermore, from (5) we can also get

$$\frac{y_1}{y_0} = -\frac{\phi_{-\gamma^-}(\tau^-)}{\phi_{\gamma^-}(\tau^-)} e^{(\gamma^- + \frac{\alpha^-}{\beta^-})\tau^-}, \tag{8}$$

where $\tau^- \in (0, \hat{\tau}^-)$ and $\hat{\tau}^-$ is the first positive solution of $\phi_{|\gamma^-|}(\hat{\tau}^-) = 0$.

Similarly, consider the system $\dot{x} = A^+ x$, take initial state $q(0, y_1, z_1)$ the expression of the solution $\Phi^+(q, t)$ can be written as

$$\Phi^+(q, t) = e^{\alpha^+ t} \cdot \tilde{M} \cdot \begin{bmatrix} \cos \beta^+ t & -\sin \beta^+ t & 0 \\ \sin \beta^+ t & \cos \beta^+ t & 0 \\ 0 & 0 & e^{(\lambda^+ - \alpha^+)t} \end{bmatrix} \cdot \tilde{M}^{-1} \cdot \begin{bmatrix} 0 \\ y_1 \\ z_1 \end{bmatrix}, \tag{9}$$

with

$$\tilde{M} = \begin{bmatrix} 1 & 0 & 1 \\ \lambda^+ + \alpha^+ & \beta^+ & 2\alpha^+ \\ \lambda^+\alpha^+ & \lambda^+\beta^+ & (\alpha^+)^2 + (\beta^+)^2 \end{bmatrix}.$$

Suppose $\Phi^+(q, \tau^+/\beta^+) = (0, y_2, z_2)$ holds for the first time with $y_1 < 0$, then by (9) we can get the corresponding parametric representations of the initial and final coordinate ratios

$$v_1(\tau^+) = \frac{z_1}{y_1} = \lambda^+ + \beta^+ \cdot [1+(\gamma^+)^2] \cdot \frac{e^{\gamma^+\tau^+}\sin\tau^+}{\phi_{\gamma^-}(\tau^-)}$$

$$v_2(\tau^+) = \frac{z_2}{y_2} = \lambda^+ - \beta^+ \cdot [1+(\gamma^+)^2] \cdot \frac{e^{-\gamma^+\tau^+}\sin\tau^+}{\phi_{-\gamma^+}(\tau^+)}$$

(10)

which not only give the angular behavior of the Poincaré half map $P^+$ but also give the definition of the slope transition map $S^+$ in the parametric form

$$S^+ : R \to R, \quad S^+(v_1) = v_2,$$

(11)

Furthermore,

$$\frac{y_2}{y_1} = -\frac{\phi_{-\gamma^+}(\tau^+)}{\phi_{\gamma^+}(\tau^+)} e^{(\gamma^+ + \frac{\alpha^+}{\beta^+})\tau^+}$$

(12)

where $\tau^+ \in (0, \hat{\tau}^+)$ and $\hat{\tau}^+$ is the first positive solution of $\phi_{|\gamma^+|}(\hat{\tau}^+) = 0$.

**Remark 2.1** *Let $\vec{n}$ be the normal vector of the separation plane $x_1 = 0$ and $\vec{g}^\pm = A^\pm x$, then by $\{(x_1, y, z) : \vec{n} \cdot \vec{g}^\pm = \delta^\pm x_1 - y = 0, x_1 = 0\} = \{(x_1, y, z) : x_1 = 0, y = 0\}$, it is not difficult to see that $S^-$ and $S^+$, and hence their composition map $S^+ \circ S^-$, are all maps having definition on the whole $R$.*

**Remark 2.2** *By $\prod_F^\pm$ denote the focus planes, i.e., the two-dimensional center invariant manifolds of the systems $\dot{x} = A^- x$ and $\dot{x} = A^+ x$, respectively. We can get*

$$\prod_F^- \equiv (\lambda^-)^2 x_1 - \lambda^- y + z = 0, \quad \prod_F^+ \equiv (\lambda^+)^2 x_1 - \lambda^+ y + z = 0.$$

(13)

*It is easy to see that if and only if $\lambda^+ = \lambda^-$ the focus plane $\prod_F^-$ coincide with the focus plane $\prod_F^+$ and $\prod_F^+$ (or $\prod_F^-$) forms a two-zonal invariant cone for system (1), which will be called the **trivial two-zonal invariant cone**. In the following analysis, without explicit directions, two-zonal invariant cones always mean the non-trivial case. Also the one-dimensional invariant manifolds of the systems $\dot{x} = A^- x$ and $\dot{x} = A^+ x$ can be given by*

$$l^- : \frac{x}{1} = \frac{y}{2\alpha^-} = \frac{z}{(\alpha^-)^2 + (\beta^-)^2}, \quad l^+ : \frac{x}{1} = \frac{y}{2\alpha^+} = \frac{z}{(\alpha^+)^2 + (\beta^+)^2}.$$

**Remark 2.3** *Note that $\tau^\pm$ (which will be called the **phase angles**) in (6)-(12) represent the*

*included angle between the projections of the initial and final half-lines passing through the origin and the points $p$ and $q$, respectively, along the one-dimensional invariant manifolds $l^{\pm}$ on the focus planes $\prod_F^{\pm}$, see **Fig. 3**.*

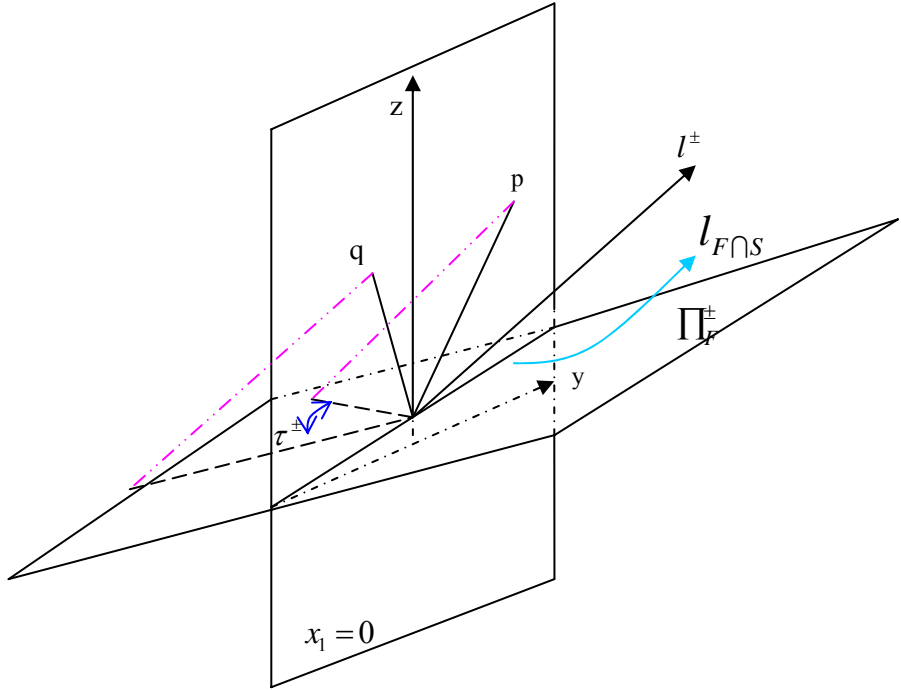

**Fig 3.** The geometric meaning of the angles $\tau^{\pm}$

## 3 Revisit to a result

The following result has been given in [7] by transferring system (1) to a continuous piecewise cubic system on $S^2$, see **Table 1**.

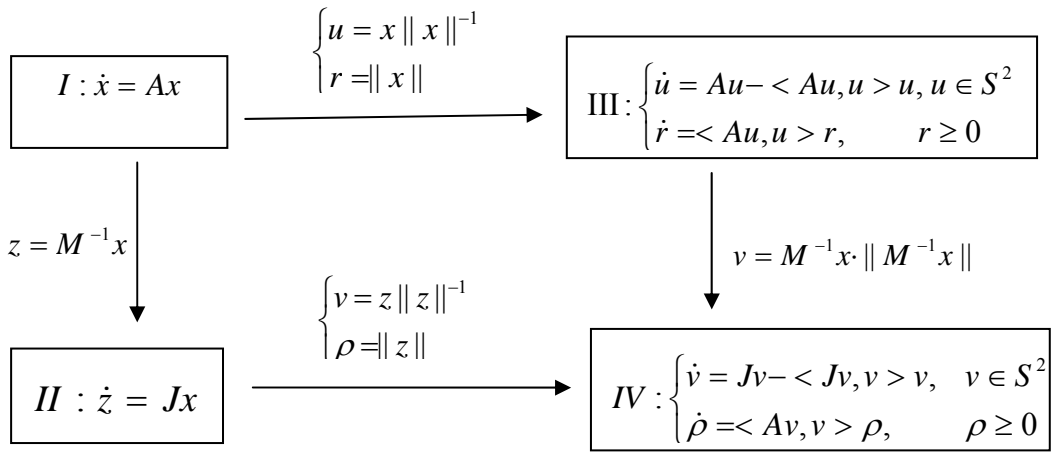

**Table 1.** The equivalent transformations between systems I -IV

Here we revisit this result to give a new proof which is more geometrically intuitive.

**Proposition 3.1 ([7])** Suppose $A$ is a $3 \times 3$ metrix given by (2), and its eigenvalues are $\lambda$, $\alpha \pm \beta j$ with $\beta > 0$, then system $\dot{x} = Ax$ has two-zonal invariant cones if and only if $\lambda = \alpha$. Furthermore, once exist, every solution starting from the separation plane $x_1 = 0$ corresponding to a two-zonal invariant cone.

**Proof** Let $\gamma = (\alpha - \lambda)/\beta$, then the solution of the system $\dot{x} = Ax$ with $(0 \; y_0 \; z_0)^T$ be the initial state is given by

$$\begin{bmatrix} x_1(t) \\ y(t) \\ z(t) \end{bmatrix} = e^{\alpha t} \cdot \overline{M} \cdot \begin{bmatrix} \cos \beta t & -\sin \beta t & 0 \\ \sin \beta t & \cos \beta t & 0 \\ 0 & 0 & e^{(\lambda - \alpha)t} \end{bmatrix} \cdot \overline{M}^{-1} \cdot \begin{bmatrix} 0 \\ y_0 \\ z_0 \end{bmatrix},$$

where

$$\overline{M} = \begin{bmatrix} 1 & 0 & 1 \\ \lambda + \alpha & \beta & 2\alpha \\ \lambda \alpha & \lambda \beta & \alpha^2 + \beta^2 \end{bmatrix}. \tag{14}$$

By $l_{F \cap S}$ denotes the intersection line between the focus plane $\Pi_F$ and the separation plane $x_1 = 0$, by (13), we know

$$l_{F \cap S} : x_1 = 0, z = \lambda y.$$

Then the proof can be finished by two cases.

**Case 1.** When $p(0 \; y_0 \; z_0)^T \notin l_{F \cap S}$. The line passing through the origin and the point $p$ can forms an invariant cone that has intersections with the plane $x_1 = 0$ under the flow of the system $\dot{x} = Ax$ if and only if its slope is the fixed point of the slope transition map $S$ defined by (7) or (11) without the subscripts, that is, the first element of the final state is zero after the time $2\pi/\beta$, see **Fig.3**, i.e.

$$x_1(2\pi/\beta) = 0,$$

which verifies

$$\frac{1}{\beta^2(1+\gamma^2)} e^{\alpha \cdot 2\pi/\beta} \cdot (1 - e^{2\pi \gamma}) = 0.$$

Hence

$$\gamma = 0 \quad \text{i.e.} \quad \lambda = \alpha.$$

**Case 2.** When $(0 \; y_0 \; z_0)^T \in l_{F \cap S}$. By the expression of the solution, we have

$$\begin{bmatrix} x_1(2\pi/\beta) \\ y(2\pi/\beta) \\ z(2\pi/\beta) \end{bmatrix} = \frac{e^{2\pi\alpha/\beta}}{\beta^3(1+\gamma^2)} \cdot \begin{bmatrix} 0 \\ 1 \\ \lambda \end{bmatrix}$$

which implies that for each initial state $(0 \ y_0 \ z_0)^T \in l_{F \cap S}$, the focus plane $\Pi_F$ itself forms a two-zonal invariant cone of the system $\dot{x} = Ax$, which will be called the **trivial two-zonal invariant cone** in our following analysis.

Because of arbitrariness of the initial state, the proof can be easily completed. □

## 4 Main results about the existence of periodic orbits

By studying the fixed points of the composition slope transition map $S^+ \circ S^-$ given by (7) and (11) in parametric forms, [6] and [7] detected the existence of the invariant cones of system (1) and have complete conclusions. Here, on the basis of this work we get some results about the existence of periodic orbits. Note that the homogeneous property of system (1) excludes the possibility of limit cycle. In fact, the dynamics on the cone is either of stable focus type, or a center, or of unstable focus type. The following theorem gives a necessary and sufficient condition for the existence of periodic orbits, that is, the dynamics on the cone is of a center type.

**Theorem 4.1** Suppose the eigenvalues of matrices $A^\pm$ are $\lambda^\pm$, $\alpha^\pm \pm \beta^\pm j$ with $\beta^\pm > 0$, then system (1) has periodic solutions if and only if there exist $\tau^\pm \in (0, \hat{\tau}^\pm)$ such that

a) $\lambda^+ - \beta^+ \cdot [1+(\gamma^+)^2] \cdot \dfrac{e^{-\gamma^+ \tau^+} \sin \tau^+}{\phi_{-\gamma^+}(\tau^+)} = \lambda^- + \beta^- \cdot [1+(\gamma^-)^2] \cdot \dfrac{e^{\gamma^- \tau^-} \sin \tau^-}{\phi_{\gamma^-}(\tau^-)}$;

b) $\lambda^+ + \beta^+ \cdot [1+(\gamma^+)^2] \cdot \dfrac{e^{\gamma^+ \tau^+} \sin \tau^+}{\phi_{\gamma^+}(\tau^+)} = \lambda^- - \beta^- \cdot [1+(\gamma^-)^2] \cdot \dfrac{e^{-\gamma^- \tau^-} \sin \tau^-}{\phi_{-\gamma^-}(\tau^-)}$;

c) $\left[ \dfrac{\phi_{-\gamma^-}(\tau^-)}{\phi_{\gamma^-}(\tau^-)} e^{2\gamma^- \tau^-} \right] \cdot \left[ \dfrac{\phi_{-\gamma^+}(\tau^+)}{\phi_{\gamma^+}(\tau^+)} e^{2\gamma^+ \tau^+} \right] \cdot e^{\frac{\lambda^+}{\beta^+}\tau^+ + \frac{\lambda^-}{\beta^-}\tau^-} = 1$,

where $\hat{\tau}^\pm$ are the first positive solutions to $\phi_{|\gamma^\pm|}(\hat{\tau}^\pm) = 0$.

**Proof** It is easily to see that system (1) has periodic orbits if and only if there exist two-zonal invariant cones and $y_0 = y_2$, which are defined in (8) and (12). Note that the sufficient and necessary condition for the existence of two-zonal invariant cones is

$$u_0(\tau^-) = v_2(\tau^+), u_1(\tau^-) = v_1(\tau^+)$$

for some $\tau^+ \in (0, \hat{\tau}^+)$ and $\tau^- \in (0, \hat{\tau}^-)$. By (6) and (10) the above equations can be rewritten as

$$\lambda^+ - \beta^+ \cdot [1+(\gamma^+)^2] \cdot \frac{e^{-\gamma^+\tau^+}\sin\tau^+}{\phi_{-\gamma^+}(\tau^+)} = \lambda^- + \beta^- \cdot [1+(\gamma^-)^2] \cdot \frac{e^{\gamma^-\tau^-}\sin\tau^-}{\phi_{\gamma^-}(\tau^-)},$$

$$\lambda^+ + \beta^+ \cdot [1+(\gamma^+)^2] \cdot \frac{e^{\gamma^+\tau^+}\sin\tau^+}{\phi_{\gamma^+}(\tau^+)} = \lambda^- - \beta^- \cdot [1+(\gamma^-)^2] \cdot \frac{e^{-\gamma^-\tau^-}\sin\tau^-}{\phi_{-\gamma^-}(\tau^-)}.$$

In addition, by (8) and (12) $y_0 = y_2$ is equivalent to

$$\left[\frac{\phi_{-\gamma^-}(\tau^-)}{\phi_{\gamma^-}(\tau^-)}e^{2\gamma^-\tau^-}\right] \cdot \left[\frac{\phi_{-\gamma^+}(\tau^+)}{\phi_{\gamma^+}(\tau^+)}e^{2\gamma^+\tau^+}\right] \cdot e^{\frac{\lambda^+}{\beta^+}\tau^+ + \frac{\lambda^-}{\beta^-}\tau^-} = 1.$$

From the above statements, the proof can be finished. □

Note that when $\gamma^+ = \gamma^- = 0$ the conditions a) and b) in **Theorem 4.1** can be reduced to

$$\lambda^+ - \beta^+ \cdot \frac{\sin\tau^+}{1-\cos\tau^+} = \lambda^- + \beta^- \cdot \frac{\sin\tau^-}{1-\cos\tau^-},$$

$$\lambda^+ + \beta^+ \cdot \frac{\sin\tau^+}{1-\cos\tau^+} = \lambda^- - \beta^- \cdot \frac{\sin\tau^-}{1-\cos\tau^-}.$$

From these two equations we can see that system (1) has two-zonal invariant cones if and only if $\lambda^+ = \lambda^-$, and every pair of phase angles $(\tau^-,\tau^+) \in \Omega$ corresponding to a two-zonal invariant cone, where

$$\Omega = \left\{(\tau^-,\tau^+) : \frac{\cot(\tau^-/2)}{\cot(\tau^+/2)} = -\frac{\beta^+}{\beta^-}\right\} \cup \{(\pi,\pi)\}.$$

Furthermore, by condition c) in **Theorem 4.1**, it is not difficult to get the following corollary.

**Corollary 4.2** Suppose $A^\pm$ in system (1) have eigenvalues $\lambda^\pm, \alpha^\pm \pm \beta^\pm j$, with $\beta^\pm > 0$, the following statements hold

a) When $\gamma^+ = \gamma^- = 0$, the dynamics on the cone is of a center (stable focus /unstable focus) type if and only if $\lambda^+ = \lambda^- = 0 (<0 />0)$;

b) The dynamics on the trivial two-zonal invariant cone is of a center (stable focus /unstable focus) type if and only if $\frac{\alpha^+}{\beta^+} + \frac{\alpha^-}{\beta^-} = 0 (<0 />0)$.

Although **Theorem 4.1** gives a necessary and sufficient condition, it is generally still very difficult to determine whether a certain system has periodic orbits or not, since it is impossible to get the analytic solution of the existence conditions given by a), b）and c) in **Theorem 4.1**. In order to illustrate how these conditions are used to show the existence of systems satisfying these conditions and hence having periodic orbits, we give the following conclusion, in which for the clarity of the discussion we fix the values of some parameters.

**Corollary 4.3** Suppose $A^{\pm}$ in system (1) have eigenvalues $\lambda^{\pm}, \alpha^{\pm} \pm \beta^{\pm} j$, with $\beta^{\pm} > 0$, Fix $\gamma^+ \cdot \gamma^- > 0$ and let $\lambda^+ - \lambda^- = c$, then for each $c \in R$ there exist a group of systems with the form given by (2) having periodic orbits.

**Proof** When $c = 0$, it is not difficult to see by **Corollary 4.2** that those systems that have conjugate complex eigenvalues $\alpha^{\pm} \pm \beta^{\pm} j$ satisfying $\dfrac{\alpha^+}{\beta^+} + \dfrac{\alpha^-}{\beta^-} = 0$ all have periodic orbits.

When $c \neq 0$, take $\gamma^+ = \gamma$ and $\gamma^- = k\gamma$, then $k > 0$. Let $\hat{\tau}^{\pm}$ be the first positive solutions for $\phi_{|\gamma^{\pm}|}(\hat{\tau}^{\pm}) = 0$. In the subsequence analysis we will show that for each $c \neq 0$ every pair of angles $(\tau^-, \tau^+) \in \widetilde{\Omega}$ determine a system with form (2) having periodic orbits, and hence the proof can be finished, where

$$\widetilde{\Omega} = \{(\tau^-, \tau^+): \tau^+ \in (0, \hat{\tau}^+), \tau^- \in (0, \hat{\tau}^-), \tau^{\pm} \neq \pi, \sin\tau^- \cdot \sin\tau^+ < 0, \mathrm{sgn}(\sin\tau^-) \cdot \mathrm{sgn}(\gamma) = \mathrm{sgn}(c)\}.$$

In order to continue, here we first introduce an auxiliary function

$$g_\gamma(\tau) = \frac{\phi_{-\gamma}(\tau)}{\phi_\gamma(\tau)} e^{2\gamma\tau}, \tau \in (0, \hat{\tau})$$

where $\hat{\tau}$ is the first positive solution for $\phi_{|\gamma|}(\hat{\tau}) = 0$. It is not difficult to show this function has the following properties by elementary calculation,

i) When $\gamma > 0$, $g'(\tau) > 0$, and $\lim\limits_{\tau \to 0^+} g(\tau) = 1$, $\lim\limits_{\tau \to \hat{\tau}^-} g(\tau) = +\infty$;

ii) When $\gamma < 0$, $g'(\tau) < 0$, and $\lim\limits_{\tau \to 0^+} g(\tau) = 1$, $\lim\limits_{\tau \to \hat{\tau}^-} g(\tau) = 0$;

iii) When $\gamma = 0$, $g(\tau) \equiv 1$ for $\tau \in (0, 2\pi)$.

The graph of $g_\gamma(\tau)$ is given in **Fig. 4**, where by the dot-dash line we denote the line $\tau = \hat{\tau}$.

Let

$$\Delta = -(1+\gamma^2) \cdot (1+k^2\gamma^2) \cdot \frac{\sin\tau^- \cdot \sin\tau^+ \cdot e^{-\gamma\tau^+} e^{k\gamma\tau^-}}{\phi_{-\gamma}(\tau^+) \cdot \phi_{k\gamma}(\tau^-)} \left[ \frac{\phi_{-\gamma}(\tau^+)}{\phi_\gamma(\tau^+)} \cdot e^{2\gamma\tau^+} - \frac{\phi_{k\gamma}(\tau^-)}{\phi_{-k\gamma}(\tau^-)} \cdot e^{-2k\gamma\tau^-} \right].$$

By the properties of the auxiliary function $g_\gamma(\tau)$ we know

$$\operatorname{sgn}(\Delta) = \operatorname{sgn}(\gamma).$$

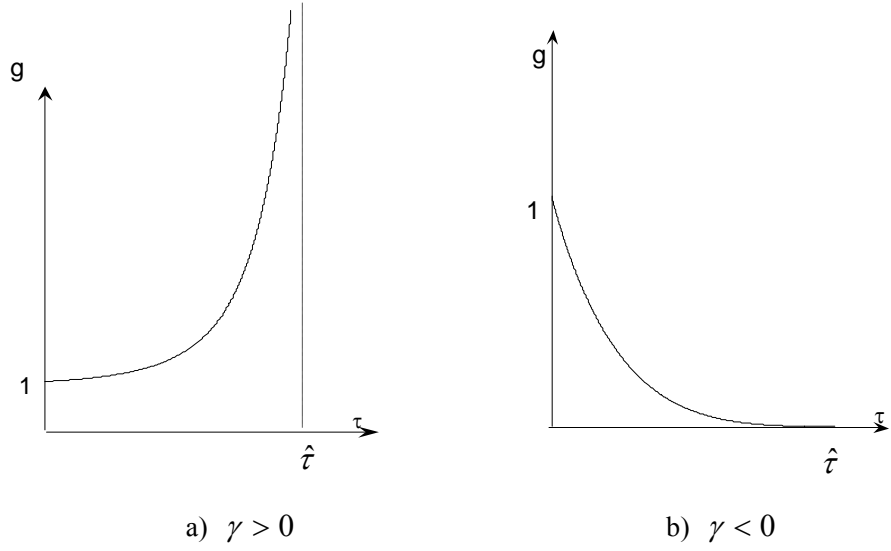

a) $\gamma > 0$      b) $\gamma < 0$

**Fig.4**   The graph of $g$ according to the sign of $\gamma$

Then by solving the set of equations a), b) and c) in **Theorem 4.1**, it is not difficult to get the values of $\beta^{\pm}$ and $\lambda^{\pm}$

$$\beta^{-} = -\frac{c}{\Delta} \cdot (1+\gamma^2) \cdot \sin \tau^{+} \cdot \left[ \frac{e^{\gamma \tau^{+}}}{\phi_{\gamma}(\tau^{+})} + \frac{e^{-\gamma \tau^{+}}}{\phi_{-\gamma}(\tau^{+})} \right] > 0,$$

$$\beta^{+} = \frac{c}{\Delta} \cdot (1+k^2\gamma^2) \cdot \sin \tau^{-} \cdot \left[ \frac{e^{k\gamma \tau^{-}}}{\phi_{k\gamma}(\tau^{-})} + \frac{e^{-k\gamma \tau^{-}}}{\phi_{-k\gamma}(\tau^{-})} \right] > 0,$$

$$\lambda^{-} = \left( \widetilde{\nabla} - \frac{\tau^{+} \cdot c}{\beta^{+}} \right) \Big/ (\frac{\tau^{+}}{\beta^{+}} + \frac{\tau^{-}}{\beta^{-}}),$$

$$\lambda^{+} = \left( \widetilde{\nabla} + \frac{\tau^{-} \cdot c}{\beta^{-}} \right) \Big/ (\frac{\tau^{+}}{\beta^{+}} + \frac{\tau^{-}}{\beta^{-}}),$$

where

$$\widetilde{\nabla} = -\ln \left( \frac{\phi_{-\gamma}(\tau^{+})}{\phi_{\gamma}(\tau^{+})} e^{2\gamma \tau^{+}} \cdot \frac{\phi_{-k\gamma}(\tau^{-})}{\phi_{k\gamma}(\tau^{-})} e^{2k\gamma \tau^{-}} \right).$$

Thus the system with $\lambda^{\pm}, \alpha^{\pm} \pm \beta^{\pm} j$ (where $\alpha^{\pm}$ can be got by substituting $\beta^{\pm}$ and $\lambda^{\pm}$ into $\gamma^{\pm} = (\alpha^{\pm} - \lambda^{\pm})/\beta^{\pm}$) be its eigenvalues is the system that has periodic orbits, and the proof is completed.□

Note that for the cases $\gamma^+ \cdot \gamma^- \leq 0$, the sign of $\nabla$ cannot be determined without additional conditions therefore the positivity of $\beta^\pm$ can neither be determined. This implies that the systems in question do not necessarily have invariant cones, and periodic orbits might not exist. Put

$$\widetilde{\widetilde{\nabla}} = \frac{\phi_{-\gamma}(\tau^+)}{\phi_\gamma(\tau^+)} e^{2\gamma\tau^+} \cdot \frac{\phi_{-k\gamma}(\tau^-)}{\phi_{k\gamma}(\tau^-)} e^{2k\gamma\tau^-}.$$

Then by virtue of the auxiliary function $g_\gamma(\tau)$ introduced in **Corollary 4.3**, it is straightforward to see that $\widetilde{\widetilde{\nabla}}$ is equal to one when $\gamma^+ = \gamma^- = 0$, greater than one when $\gamma^+ \geq 0, \gamma^- \geq 0$ with $\gamma^+ + \gamma^- \neq 0$, and less than one when $\gamma^+ \leq 0, \gamma^- \leq 0$ with $\gamma^+ + \gamma^- \neq 0$.

By condition c) of **Theorem 4.1**, we have

**Corollary 4.4.** Suppose the eigenvalues of matrices $A^\pm$ in system (1) are $\lambda^\pm$, $\alpha^\pm \pm \beta^\pm j$ with $\beta^\pm > 0$, then we have the following statements

a) When $\gamma^+ = \gamma^- = 0$, then the necessary condition for the existence of periodic orbits is $\lambda^+ \cdot \lambda^- \leq 0$;

b) When $\gamma^+ \geq 0, \gamma^- \geq 0$ and $\gamma^+ + \gamma^- \neq 0$, then the necessary condition for the existence of periodic orbits is that at least one of $\lambda^+, \lambda^-$ is negative;

c) When $\gamma^+ \leq 0, \gamma^- \leq 0$ and $\gamma^+ + \gamma^- \neq 0$, then the necessary condition for the existence of periodic orbits is that at least one of $\lambda^+, \lambda^-$ is positive;

## 4 Examples

**Example 1.** In order to illustrate **Theorem 4.1**, we consider (1) with specific parameters. Let $\gamma^+ = 1$, $\gamma^- = 1$, then the first positive solutions to $\phi_{|\gamma^\pm|}(\hat{\tau}^\pm) = 0$ satisfies $\hat{\tau}^\pm \in (\frac{5\pi}{4}, 2\pi)$. Take $\tau^+ = \frac{5\pi}{4}, \tau^- = \frac{\pi}{4}$ and $c = 10$, by virtue of the proof of **Corollary 4.3**, we get

$$A^- = \begin{pmatrix} -24.5442 & -1 & 0 \\ 207.7430 & 0 & -1 \\ -629.2483 & 0 & 0 \end{pmatrix}, \quad A^+ = \begin{pmatrix} -0.5542 & -1 & 0 \\ 0.1349 & 0 & -1 \\ -0.0203 & 0 & 0 \end{pmatrix}$$

with eigenvalues

$$\lambda^- = -10.3322, \quad \alpha^- \pm \beta^- j = -7.1060 \pm 3.2259 j,$$

$$\lambda^+ = -0.3321, \quad \alpha^+ \pm \beta^+ j = -0.1111 \pm 0.2209 j.$$

It is easy to verify that $\gamma^\pm, \tau^\pm, \lambda^\pm$ and $\beta^\pm$ given above satisfy the conditions a),b) and c) in **Theorem 4.1**, hence the corresponding system has periodic orbits, as shown in **Fig.5**, where the periodic orbit with initial condition $x_0 = [0 \ 4 \ -1.3040]^T$ intersects the separation plane $x = 0$ at the points $A, B$. The segments of the periodic orbit with $x > 0, x < 0$ are represented by dot line and solid line, respectively.

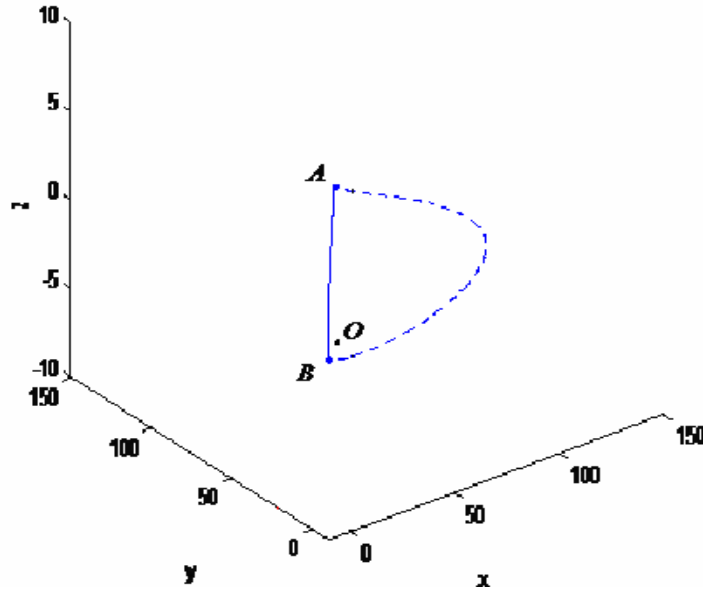

**Fig.5** A periodic orbit of system (1) with $x_0 = [0 \ 4 \ -1.3040]^T$

**Example 2.** Take $\gamma^+ = 1, \gamma^- = 1, c = -10$ and $\tau^- = \dfrac{5\pi}{4}, \tau^+ = \dfrac{\pi}{4}$, then

$$A^- = \begin{pmatrix} -0.5542 & -1 & 0 \\ 0.1349 & 0 & -1 \\ -0.0203 & 0 & 0 \end{pmatrix}, \quad A^+ = \begin{pmatrix} -24.5442 & -1 & 0 \\ 207.7430 & 0 & -1 \\ -629.2483 & 0 & 0 \end{pmatrix}$$

whose eigenvalues are

$$\lambda^+ = -10.3322, \quad \alpha^+ \pm \beta^+ j = -7.1060 \pm 3.2259 j.$$

$$\lambda^- = -0.3321, \quad \alpha^- \pm \beta^- j = -0.1111 \pm 0.2209 j.$$

Similarly, $\gamma^\pm, \tau^\pm, \lambda^\pm$ and $\beta^\pm$ here satisfy the conditions a),b) and c) in **Theorem 4.1**, and the periodic orbit of system (1) with initial condition $x_0 = [0 \ -4 \ 1.3040]^T$ is given in **Fig.6**.

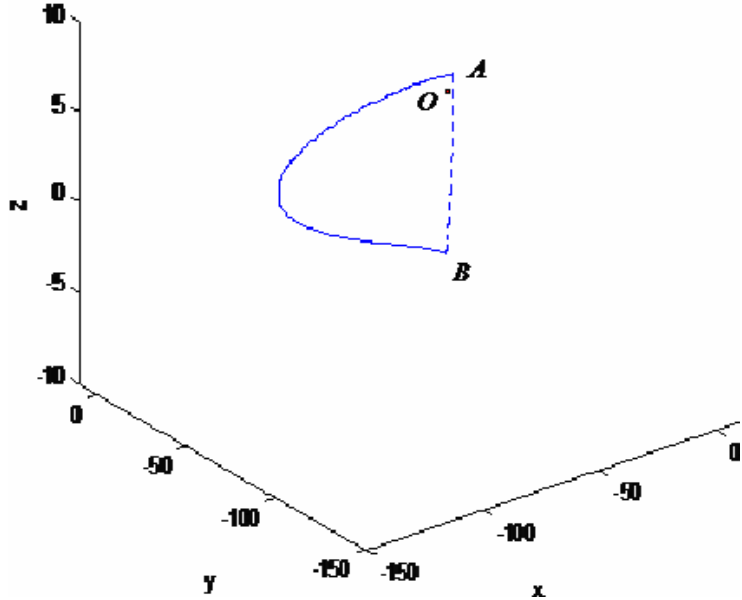

**Fig.6**  A periodic orbit of system (1) with $x_0 = [0 \ -4 \ 1.3040]^T$

**Acknowledgements** This work is supported in part by National Natural Science Foundation of China (10972082).

**References**

1 A.A. Andronov, A. Vitt and S. Khaikin, *Theroy of Oscillations*, Pergamon Press, Oxford, 1966.

2 A. Arneodo, P. Coullet and C. Tresser, *Possible new strange attractors with spiral structure.* Comun. Math. Phys [1981]. 79, 573-579.

3 B. Brogliato, *Impacts in mechanical systems-Analysis and Modeling.* Springer-Verlag, New York, Lecture notes in physics, Vol.551 (2000).

4 B. Brogliato, *Nonsmooth mechanics-models, dynamics and control.* Springer-Verlag, New York, (1999).

5 M. K. Camlibel, W. P. M. H. Heemels and J.M. Schumacher, *Stability and controllability of planar bimodal linear complementarity systems*, Proceedings of the 42$^{nd}$ IEEE Conference on Decision and Control Maul, Hawall USA(2003), 1651-1656.

6 V.Carmona, E.Freire, E.Ponce and F.Torres, *The Continuous Matching of Two Stable Linear Systems Can be Unstable*. Discrete and Continuous Dynamical Systems, Vol. 16, No.3 (2006), 689-703.

7 V.Carmona, E.Freire, E.Ponce and F.Torres, *Bifurcation of Invariant Cones in Piecewise Linear Homogeneous Systems*. International Journal of Bifurcation and Chaos, Vol. 15, No.8 (2005), 2469-2484.

8 V.Carmona, E.Freire, E.Ponce and F.Torres, *Invariant manifolds of periodic orbits for piecewise linear three-dimensional systems*. IMA J. APPl.Math, 69(2004), 71-91.

9 V.Carmona, E.Freire, E.Ponce and F.Torres, *On Simplifying and Classifying Piecewise-Linear Systems*. IEEE Trans. Circuites Syst.- I: Fund. Th. Appl. Vol.49, No.5 (2002), 609-620.

10 M. di Bernardo, C.J. Budd, A.R. Champneys and P. Kowalcayk, *Piecewise-smooth dynamical systems*, Theory and Application, Springer-Verlag, London (2008).

11 M di Bernardo, C.J. Budd, A.R. Champneys, P. Kowalcayk, A.B. Nordmark, G.O. Tost and P.T. Piiroinen, *Bifurcations in Nonsmooth Dynamical Systems*, SIAM Review, Vol.50 (2008), 629-701.

12 E. Freire, E. Ponce and J. Ros, *Limit cycle bifurcation from center in symmetric piecewise linear systems*. International Journal of Bifurcation and Chaos. 9 (1999), 895-907.

13 E. Freire, E. Ponce, E. Rodrigo and F. Torres, *Bifurcation sets of continuous piecewise linear systems with two zones*. International Journal of Bifurcation and Chaos. Sci Engrg. 8(1998), 2073-2097.

14 E. Freire, E. Ponce, E. Rodrigo and F. Torres, *Bifurcation sets of symmetrical continuous piecewise linear systems with two zones*. International Journal of Bifurcation and Chaos. Sci Engrg. 12(2002),1675-1702.

15 G. Kriegsmann, *The rapid bifurcation of the Wien bridge oscillator, IEEE Transactions Circuits and Systems,* CAS, 34(1987), 1093- 1096.

16 M. Kunze, *Non-smooth dynamical systems*. Springer-Verlag, Berlin, Lecture notes in Mathematics, Vol.1744 (2000).

17 J. Llibre and E. Ponce, *Phase portraits of planar control systems*, Nonlinear Analysis (1996), 1177-1197.